\documentclass{article}
\usepackage{amsthm,amsmath,amsfonts,latexsym,amssymb}
\usepackage{mathrsfs,dutchcal,arydshln,gensymb,appendix,abstract}
\usepackage{geometry}
\usepackage{enumitem}
\usepackage{lipsum}
\usepackage{graphicx} 
\allowdisplaybreaks[4]
\begin{document}

\title{Existence of multiple constant mean curvature hypersurfaces for varying Riemannian metrics}

\author{XIAOXIANG JIAO AND WENDUO ZOU}
\date{}

\maketitle
\begin{abstract}
Given a closed Riemannian manifold  $(M^{n+1},g)$,$3\leq n+1\leq7$.In this paper,we will prove that for any $c>0$,suppose the number of closed $c-CMC$ hypersurfaces is finite,then there exists a metric $h$ on $M$ such that the $c-CMC$ hypersurfaces in $(M,g)$ are also $c-CMC$ hypersurfaces in $(M,h)$ and the number of $c-CMC$ hypersurfaces in $(M,h)$ is strictly greater than the number of $c-CMC$ hypersurfaces in $(M,g)$.Moreover,we will give a precise upper bound for the $L^{\frac{n+1}{2}}$ norm of $(g-h)$,which depends on the metric $g$ and the number of  $c-CMC$ hypersurfaces in $(M,g)$.
\end{abstract}

\section{Introduction}

Let $M$ be a Riemannian manifold,a hypersurface $\sum$ in $M$ is called  a $c$-constant mean curvature hypersurface ($c-CMC$ hypersurface) if it is a critical point of the functional:
\begin{equation}
\mathcal{A}^{c}=Area-cVol.\tag{1.1}
\end{equation}
If $c=0$,such a hypersurface is called a minimal hypersurface which is a critical point of the area functional.

\vspace{16pt} 
The existence of minimal hypersurfaces and CMC hypersuraces is a fundamental problem in differential geometry.Min-max theory for minimal hypersurfaces is a powerful tool for showing the existence of closed minimal hypersurfaces.The works of Almgren[1],pitts[15],and Schoen-Simon [16] have shown that every closed Riemannian manifold $M^{n+1}$,$3\leq n+1\leq7$,contains a smooth,closed minimal hypersurfaces.
In[26],Zhou and Zhu developed the min-max theory for CMC hypersurfaces and constructed closed CMC hypersurfaces in closed Riemannian manifold of any prescirbed constant mean curvature.

\vspace{16pt} 
Another important problem is  the existence of multiple minimal hypersurfaces and CMC hypersurfaces in closed manifold.In [20],Yau conjectured that every closed Riemannian manifold contains infinitely many closed minimal surfaces.By combining the min-max theory of minimal hypersurfaces and  volume spectrum of Riemannian manifold which was introduced by Gromov[5,6]and was proved by Liokumovich-Marques-Neves[9] that it satisfies a weyl law,Irie-Marques-Neves[7]showed that Yau$^{,}$s conjecture is true for generic metric on closed manifold $M^{n+1} (3\leq n+1\leq7)$.Later Song[18] proved that  Yau$^{,}$s conjecture is true for any metric on closed manifold $M^{n+1} (3\leq n+1\leq7)$.In higher dimension case,Li[8] proved that any closed manifold $M^{n+1}(n+1\geq8)$ contains infinitely many closed minimal hypersurfaces with the Hausdorff dimension of the singular set is less than or equal to $n-7$.

\vspace{16pt} 
In the case of CMC hypersurfaces,Dey[3] made use of the gap between the $k$-width and $(k+1)$-width to construct some homotopy classes that satisfy the homotopy conditions of min-max theory of CMC hypersurfaces with small mean curvature,he showed that the lower bound of the number of $c-CMC$ hypersurfaces increases as the mean curvature tends to zero.In the large mean curvature regime,Ye[21] proved that $M$ contains at least $\mu(M)$ $c-CMC$ spheres if $c$ is large enough and the scalar curvature function of $M$ is Morse.Where  $\mu(M)$ is the Morse number of $M$.Later,Pacard and Xu[14] proved that  $M$ contains at least cat($M$) $c-CMC$ spheres when c is large enough without the hypothesis of scalar curvature function above.Here cat($M$) is the Lusternik-Schnirelmann category of $M$.However,by the results of Mazzeo and Pacard[13],those CMC spheres do not exhaust the CMC hypersurfaces with large mean curvature.

\vspace{16pt}
Recently Mazurowski and Zhou[10] introduced the half-volume spectrum of Riemannian manifold and proved that half-volume spectrum also satisfies a Wely law,combined with the min-max theory of a novel functional developed by them in [11],they showed that for a generic metric on a closed manifold $M^{n+1} (3\leq n+1\leq5)$,there exist infinitely many distinct closed CMC hypersurfaces in $M$,and each encloses half of the volume of $M$.Furthermore,Mazurowski and Zhou [12] proved that any closed manifold $M^{n+1}$ ($3 \leq n+1 \leq 7$) equipped with a generic metric satisfies the following: for any constant $c > 0$, $M^{n+1}$ either contains infinitely many $c$-CMC hypersurfaces, or contains infinitely many CMC hypersurfaces with mean curvature less than $c$ but enclosing half the volume of $M$.

\vspace{16pt} 
However, since the min-max theory of CMC hypersurfaces is operated on the space of Caccioppoli sets,which is a contractible space.This leads to the homotopy conditions in the CMC min-max theory being much harder to achieve than the homotopy conditions in the min-max theory of minimal hypersurfaces.Thus the existence of multiple CMC hypersurfaces in closed manifold is hard to handle.
It still left open the question of the how many Riemannian metrics on a closed manifold admits multiple  $c-CMC$ hypersurfaces for any given constant $c$.It was conjectured by Zhou[25] that every closed Riemannian manifold $M^{n+1} (3\leq n+1\leq7)$ admits at least two distinct closed hypersurfaces with constant mean curvature $c$ for every $c>0$.This can be seen as a higher dimension analog of $\text{Arnold}^{,}$s conjecture[2] of existence of closed curves with constant geodesic curvature in Riemannian $2$-sphere.

\vspace{16pt} 

In this paper,we show that for any Riemannian metrics $g$ on a closed smooth manifold $M^{n+1}(3\leq n+1 \leq 7)$ and  any constant $c>0$,suppose the number of closed $c-CMC$ hypersurfaces is finite,then we can perturb $g$ to get a metric $h$ such that $(M,h)$ contains more $c-CMC$ hypersurfaces than $(M,g)$.Moreover,we will give a precise upper bound for the $L^{\frac{n+1}{2}}$ norm of $(g-h)$,which depends on the metric $g$ and the number of  $c-CMC$ hypersurfaces in $(M,g)$.

 \vspace{16pt}
 We now state the main theorem and some corollaries.

 \vspace{1.5pt}
 Let $(M^{n+1},g)$ be a closed Riemannian manifold with $3\leq n+1 \leq 7$.We use the notation $\mathcal{C}(M)$ to indicate the space of Caccioppoli sets.Let \text{\small$\Pi$} be the class of all $\mathcal{F}$-continuous maps $\phi:[0,1]\to \mathcal{C}(M)$ satisfy 
 \begin{enumerate}
     \item [1)] $\phi(0)=0,\phi(1)=M$;
    \item[\textbf{2)}] $\phi$ has no cencentration of mass,i.e
\begin{equation*}
    \lim_{r \to 0}\sup \{||\partial \phi(x)||(B_{r}(p))\:| \:  p \in M,x\in [0,1]\}=0. 
\end{equation*}
  
 \end{enumerate}
 The one parameter Almgren-Pitts width of  $(M,g)$ is defined to be:
 \begin{equation*}
   W(M,g)= \underset{\phi\ \in \text{\small$\Pi$}}{\inf} 
    \sup_{x \in [0,1]} \mathbb{M}(\partial\phi(x)).  
 \end{equation*}
 and the $\mathcal{A}^{c}$-min-max value of \text{\small$\Pi$} is defined to be:
 \begin{equation*}
   W^{c}(M,g)= \underset{\phi\ \in \text{\small$\Pi$}}{\inf} 
    \sup_{x \in [0,1]} [\mathbb{M}(\partial\phi(x))-c\mathcal{H}^{n+1}(\phi(x))].    
 \end{equation*}

 \vspace{5pt}
For any $c>0$,we denote by $\mathcal{P}^{c}_{g}$ the class of smooth,closed,almost embedded $c-CMC$ hypersurface in $(M,g)$ and $\widetilde{\mathcal{P}}^{c}_{g}$ represents the class of open set $\Omega$ such that $\partial \Omega$ is a smooth,closed,almost embedded $c-CMC$ hypersurface with respect to the inward unit normal of $\Omega$ in $(M,g)$.

\vspace{16pt}
 \textbf{Theorem 1.1}(Main Theorem).Let $(M^{n+1},g)$ be a closed Riemannian manifold of dimension $3\leq n+1 \leq7$.Suppose $\#\mathcal{P}^{c}_{g}$ is finite,and $\#\ \widetilde{\mathcal{P}}^{c}_{g}=N$.Then for any $k\in \mathbb{N}^{+}$,there exists a Riemannian metric $h$ on $M$ and a constant $a(n)$ depending only the dimension of $M$ such that:
 \begin{enumerate}
    \item[\textbf{(1)}] $\#$ $\widetilde{\mathcal{P}}^{c}_{h}\geq N+k$;
\item[\textbf{(2)}] $\mathcal{P}^{c}_{g}\subset$ $\mathcal{P}^{c}_{h}$;
\item[\textbf{(3)}] $||h-g||_{L^{\frac{n+1}{2}}(M,g)}\leq  a(n)[\frac{W^{c}(M,g)}{c}]^{\frac{2}{n+1}}k(N+\frac{k-1}{2})^{\frac{2}{n+1}} $.
 \end{enumerate}
 
 \vspace{16pt} 
 \textbf{1.2.Remark}:Proposition 3.4 in this paper states that  $W^{c}(M,g)$ depends continuously on the metric $g$.Furthermore,the magnitude of the perturbation form $g$ to $h$ increases as the number of CMC hypersurfaces increases, and decreases as the mean curvature $c$ increases.

 \vspace{16pt}
Since $W(M,g)$ coincides with $\omega_{1}(M,g)$ when  $(M,g)$ has positive Ricci curvature [22,23],where $\omega_{1}(M,g)$ is the $1-width$ of $(M,g)$.Then by  Corollary 1.3 in [19] (see also Theorem 1.3 in [4]) we have the following corollary:

\vspace{16pt}
\textbf{Corollary 1.3}.With the hypothesis and notations in Theorem 1.2,suppose $(M,g)$ has positive Ricci curvature,then there exists a constant $a(n)$ depending only the dimension of $M$ such that $||h-g||_{L^{\frac{n+1}{2}}(M,g)}\leq a(n)[\frac{(vol(M,g))^{\frac{n}{n+1}}}{c}]^{\frac{2}{n+1}}k(N+\frac{k-1}{2})^{\frac{2}{n+1}} $. 

\vspace{16pt}
By Remark 1.2 we have the following corollary holds.

\vspace{16pt}
\textbf{Corollary 1.4}.Let $(M^{n+1},g)$ be a Riemannian manifold with dimension $3\leq n+1 \leq 7$,then there exist a sequence of Riemannian metric $\{g_{i}\}_{i\in \mathbb{N}}$ and two sequences of positive real number $\{c_{i}\}_{i\in \mathbb{N}}$,$\{N_{i}\}_{i\in \mathbb{N}}$ with both two sequences tend to infinity such that 
\begin{enumerate}
    \item [1)] $g_{i}$ converges to $g$ in $L^{\frac{n+1}{2}}$-topology;
    \item [2)] $\#$ $\widetilde{\mathcal{P}}^{c_{i}}_{g_{i}}\geq N_{i}$.
\end{enumerate}

\vspace{16pt}
\textbf{1.5 Outline of Paper}. In Section 2,we first recall some basic notions used in this paper.We then review the min-max theory of CMC  hypersurfaces developed by Zhou-Zhu[26].

\vspace{5pt}
In Section 3,we give the proof of Theorem 1.1.The proof is divided into three parts.In the first part,we show that the  the $\mathcal{A}^{c}$-min-max value $W^{c}(M,g)$ depends continuously on the Riemannian metric.In the second part,we show that Theorem 1.1 holds when  $k=1$ via a suitable conformal  perturbation of the metric.In the third part,we construct a sequence of metrics $\{g_{i}\}_{i=1,...,m(k)}$ such that for each $g_{i}$,the number of $c-CMC$ hypersurfaces in $(M,g_{i})$ is strictly greater than the number  of $c-CMC$ hypersurfaces in $(M,g_{i-1})$,and  $||g_{i}-g_{i-1}||_{L^{\frac{n+1}{2}}(M,g_{i})}$ admits a suitable upper bound, we further prove that $||g_{i}-g_{i-1}||_{L^{\frac{n+1}{2}}(M,g_{i})}$ concides with  $||g_{i}-g_{i-1}||_{L^{\frac{n+1}{2}}(M,g)}$.Finally using the subadditivity of the norm,we deduce that Theorem 1.1 holds.

\vspace{20pt} 
\textbf{1.6 Acknowledgements}.This work is Supported by the National Natural Science Foundation of China(Grant no. 12371055).

\vspace{30pt}
\section{Notation and Preliminaries}

\vspace{25pt} 
{\large \textbf{2.1 Notation}} 

\vspace{10pt} 
Let $M^{n+1}$ be a closed Riemannian manifold of dimension $3\leq n+1\leq7$.We first collect some notions in geometric measure theory.For interested readers,we refer to [17] and [15] for further materials.

\vspace{10pt}
\textbf{a)} $\mathcal{H}^{n}$: $n$-dimensional Hausdorff measure;

\vspace{2pt}
\textbf{b)} $\mathcal{C}(M)$: the space of Caccioppoli sets in $M$;

\vspace{2pt}
\textbf{c)} $\mathcal{Z}_{n}(M,\mathbb{Z}_{2})$: the space of $n-$dimensional mod 2 flat cycles in $M$.

\vspace{2pt}
\textbf{d)}  $\mathcal{V}_{n}(M)$: the space of $n$-dimensional varifolds on $M$.

\vspace{2pt}
\textbf{e)} $|T|$: the varifold induced by $T\in \mathcal{Z}_{n}(M,\mathbb{Z}_{2})$.

\vspace{2pt}
\textbf{f)} $\partial \Omega$: the (reduced)-boundary of Caccioppoli set $\Omega$ as a current.

\vspace{2pt}
\textbf{g)} $\mathbb{M}_U$: the mass norm in an open set $U$ on the space of flat cycles.($\mathbb{M}$ deontes the $\mathbb{M}_M)$.

\vspace{2pt}
\textbf{h)} $\mathcal{F}$: flat norm on $\mathcal{Z}_{n}(M,\mathbb{Z}_{2})$ and $\mathcal{C}(M)$.

\vspace{2pt}
\textbf{i)} $\mathbf{F}$: the varifold $\mathbf{F}$-metric on $\mathcal{V}_{n}$$(M)$ and $\mathbf{F}$-metric on  $\mathcal{C}(M)$.By convention the  $\mathbf{F}$-metric on $\mathcal{C}(M)$ is
\begin{equation*}
\mathbf{F}(\Omega_{1},\Omega_{2})=\mathcal{F}(\Omega_{1},\Omega_{2})+\mathbf{F}(|\partial \Omega_{1}|,|\partial \Omega_{2}|)
\end{equation*}

\vspace{2pt}

\textbf{j)} $||V||,||T||$:the associated Radon measure on $M$ of $V\in \mathcal{V}_{n}(M)$ and $T\in \mathcal{Z}_{n}(M,\mathbb{Z}_{2})$.

\vspace{2pt}

\textbf{k)} $B_{r}(p):$ the geodesic ball of $(M,g)$.

\vspace{15pt} 
{\large \textbf{2.2 Almost embedded hypersurface and $\mathcal{A}^{c}$functional }} 

\vspace{15pt} 
 \textbf{Definition 2.1} (Almost embedded hypersurface). Let $\phi:\scalebox{0.9}{\(\sum\)}^{n} \rightarrow M$ be a smooth immersion,we say that $\phi:\scalebox{0.9}{\(\sum\)}^{n}$ is an almost embedding hypersurface if for any $p$$\in$ $\phi(\scalebox{0.9}{\(\sum\)})$ where $\phi$ fails to be embedded,there exists a neighborhood $U$ of $p$ such that:
 
\vspace{5pt} 
 \textbf{1)} $\scalebox{0.9}{\(\sum\)} \cap \phi^{-1}(U)$ is  decomposed into a disjoint union of finite pieces $\scalebox{0.9}{\(\sum\)}_{i},i=1,...l$,and $\phi$ restrect on each $\scalebox{0.9}{\(\sum\)}_{i}$ is an embedding;
  
  \vspace{5pt} 
  \textbf{2)} for each $i$,any other $\phi(\scalebox{0.9}{\(\sum\)}_{j})$,$j\neq i$,lies on one side of  $\phi(\scalebox{0.9}{\(\sum\)}_{i})$;

\vspace{12pt} 
Remark 2.2. By the maximum principle of CMC hypersurface in[4],we know that an almost embedded CMC hypersurfaces locally only decomposes into two embedded pieces.Moreover,if we denote the points in $\phi(\scalebox{0.9}{\(\sum\)})$ where fails to be embedded as $\mathcal{S(\scalebox{0.9}{\(\sum\))}}$,then $\mathcal{S(\scalebox{0.9}{\(\sum\))}}$ is contained in a countable union of connected embedded $(n-1)$-dimensional submanifolds.In particular,$\phi(\scalebox{0.9}{\(\sum\)})\setminus \mathcal{S(\scalebox{0.9}{\(\sum\))}}$ is open and dense in $\phi(\scalebox{0.9}{\(\sum\)})$.

  \vspace{12pt}
Given $c>0$,define the $\mathcal{A}^{c}$-functional on  $\mathcal{C}(M)$ as 
\begin{equation*}
\mathcal{A}^{c}(\Omega)=\mathcal{H}^n(\partial \Omega)-\mathcal{H}^{n+1}(\Omega).\tag{2.1}
\end{equation*}
The first variation for $\mathcal{A}^{c}$ along vector field $X$ is given by
\begin{equation}
\delta \mathcal{A}^{c}|_{\Omega}(X) = \int_{\partial \Omega} \text{div}_{\partial \Omega} X\, d\mathcal{H}^n - c \int_{\partial \Omega} g(X,v_{\partial \Omega}) \, d\mathcal{H}^n
\tag{2.2}
\end{equation}
where $v_{\partial \Omega}$ is the outward normal on $\partial \Omega$.Suppose $\scalebox{0.9}{\(\sum\)}$=$\partial \Omega$ is a smooth hypersurface of $M$ and $\Omega$ is a critical point of $\mathcal{A}^{c}$,then by the first variation formula of volume,we have $\scalebox{0.9}{\(\sum\)}$ has constant mean curvature c with respect to the $-v_{\partial \Omega}$.

\vspace{15pt} 
{\large \textbf{2.3 Min-Max theory for CMC hypersurfaces}} 

\vspace{15pt} 
We first recall some basic notations in min-max theory of CMC hypersurfaces.

\vspace{12pt}
\textbf{Definition 2.3}.Let $\phi:$$[0,1] \to \mathcal{C}(M)$ be a 
$\mathcal{F}$-continuous map,we call $\phi$ a \textit{sweepout} if $\phi$ satisfies the following conditions:
\begin{enumerate}
    \item[\textbf{1)}] $\phi(0)=0,\phi(1)=M$;
    \item[\textbf{2)}] $\phi$ has no cencentration of mass,i.e
\begin{equation}
      \lim_{r \to 0}\sup \{||\partial \phi(x)||(B_{r}(p))\:| \:  p \in M,x\in [0,1]\}=0.  \tag{2.4}
\end{equation}
\end{enumerate}

\vspace{9pt}
Let \text{\small$\Pi$} denotes the  collection of all sweepouts of $M$.Now we can define the $\mathcal{A}^{c}$-min-max value of \text{\small$\Pi$}.

\vspace{12pt}
\textbf{Definition 2.4}.For every $c>0$,the $\mathcal{A}^{c}\text{-}min\text{-}max\ value$ of \text{\small$\Pi$}  is defined by:
\begin{equation*}
    W^{c}(M,g)= \underset{\phi\ \in \text{\small$\Pi$}}{\inf} 
    \sup_{x \in [0,1]} \mathcal{A}^{c}(\phi(x)).
\end{equation*}

\vspace{12pt}
\textbf{Proposition 2.5}.For any metric $g$ on $M$,we have $W^{c}(M,g)$ is positive.
\begin{proof}
By the interpolation results of Zhou[23,Theorem 5.1](see also [24,Theorem 1.11 and Theorem 1.12]),we can further require that sweepout is continuous in $\mathbf{F}$-topology.Then the proposition follows from the isopermetric inequality as in the Remark 3.10 in [26].
\end{proof}

\vspace{12pt}
\textbf{Definition 2.6}.A sequence $\{\phi_{i}\}_{i\in \mathbb{N}}\in$\text{\small$\Pi$} is called a $min\text{-}max\ sequence$ if
\begin{equation*}
   \limsup_{i\to \infty}[\sup_{x\in[0,1]}\mathcal{A}^{c}(\phi_{i}(x))]=W^{c}(M,g)  
\end{equation*}
If $\{\phi_{i}\}_{i\in \mathbb{N}}$ is a min-max sequence,the $critical\ set$ of $\{\phi_{i}\}$ is defined by
\begin{equation*}
    \textbf{C}(\{\phi_{i}\})=\{V \in \mathcal{V}_{n}(M)|V=\lim_{j\to \infty}|\partial \phi_{i_{j}}(x_{j})|\
     \ \text{as varifolds with $\lim_{j\to \infty }\mathcal{A}^{c}( \phi_{i_{j}}(x_{j}))=W^{c}(M,g)$}\} 
\end{equation*}

\vspace{12pt}
The following min-max theorem essentially follows from [26,Theorem 3.8].The discrete version of min-max theory of CMC hypersurfaces was firstly proved by Zhou and Zhu in[26],and the continuous version  in higher dimension was proved by Dey in [3].We will only outline the necessary changes to the proof of Theorem 3.8 in [26].

\vspace{12pt}
\textbf{Theorem 2.7}([26,Theorem 3.8]).\textit{Let $M^{n+1}(3 \leq n+1 \leq 7)$ be a closed Riemannian manifold and $c>0$.Then there exists an open set  $\Omega$ such that $\scalebox{0.9}{\(\sum\)}$=$\partial \Omega$ is a nontrivial,smooth,closed,almost embedded c-CMC(with respect to the inward unit normal) hypersurface.Moreover, $\mathcal{A}^{c}(\Omega) = W^{c}(M,g)$.}
\begin{proof}
By the interpolation results of Zhou[23,Theorem 5.1](see also [24,Theorem 1.11 and Theorem 1.12]),there exists a min-max sequence $\{\phi_{i}\}_{i \in \mathbb{N}} $ continuous in the $\mathbf{F}$-topology.Then the pulltight procedure  [26,Proposition 4.4] works for $\{\phi_{i}\}_{i \in \mathbb{N}} $.Thus we obtain a min-max sequence $\{\Phi_{i}\}_{i \in \mathbb{N}} $ continuous in the $\mathbf{F}$-topology and each $V\in \textbf{C}(\{\Phi_{i}\})$ has $c$-bounded first variation.Since the deformation process in [26,Theorem 5.6] was only made to the points $\Phi_{i}(x)$ close to $W^{c}(M,g)$,then by Proposition 2.5 we deduce that after the deformation process in [26,Theorem 5.6] we still get a sequence of sweepout.Then we complete the proof by applying the Theorem 6.1 in [26].
\end{proof}
\section{Proof of the main theorems}

\vspace{15pt}
Let $\Omega$ be a Caccioppoli set,let $\mathbb{M}_U(\partial \Omega)_{g}$ denotes the mass of $\partial \Omega$ in an open subset $U\subset M$ with respect a Riemannian metric $g$.($\mathbb{M}(\partial \Omega)_{g_{1}}$ denotes the $\mathbb{M}_M(\partial \Omega)_{g}$).We first prove that the definition of Caccioppoli set is independent of choice of metrics.

\vspace{12pt}  
\textbf{Lemma 3.1}. \textit{Let $g_{1},g_{2}$ be two Riemannian metrics on $M$. Suppose $\Omega$ is a Caccioppoli set in $(M,g_{1})$, then $\Omega$ is also a Caccioppoli set in $(M,g_{2})$.}
\begin{proof}
Suppose $\Omega$ is a Caccioppoli set in $(M,g_{1})$, let $X$ be a vector field such that $\|X\|_{g_{2}} \leq 1$.Let $\mathcal{H}^{n+1}_{g_{i}}$ denote the $(n+1)-$dimensional Hausdorff measure with respect to ${g_{i}}$,since for every measurable subset $A$,we have 
\begin{equation}
\mathcal{H}^{n+1}_{g_{i}}(A)=\int_{A } \frac{\sqrt{\det(g_{i})}}{\sqrt{\det(g_{j})}}d\mathcal{H}^{n+1}_{g_{j}}.\tag{3.1}
\end{equation}
Where $\frac{\sqrt{\det(g_{i})}}{\sqrt{\det(g_{j})}}$ can be represented as $\frac{\sqrt{\det(g_{i}(\partial x_{k},\partial x_{l})}}{\sqrt{\det(g_{j}(\partial x_{k},\partial x_{l})}}$ in local coordinate on $M$ .Note that $\frac{\sqrt{\det(g_{i})}}{\sqrt{\det(g_{j})}}$ is a well defined continuous function on $M$, hence $ \sup_{x \in M}\frac{\sqrt{\det(g_{i})(x)}}{\sqrt{\det(g_{j})(x)}} $ is finite,thus we have
\begin{align}
\int_{\Omega } \text{div}_{g_{2}} X \, d\mathcal{H}^{n+1}_{g_{2}} 
&= \int_{\Omega } \text{div}_{g_{1}} \left( \frac{\sqrt{\det(g_{2})}}{\sqrt{\det(g_{1})}} X \right) d\mathcal{H}^{n+1}_{g_{1}} \notag \\
&\leq \sup_{x \in M}\frac{\sqrt{\det(g_{2})(x)}}{\sqrt{\det(g_{1})(x)}} \mathbb{M}(\partial \Omega)_{g_{1}}. \tag{3.2}
\end{align}
\end{proof}
By the proof above,we have the following corollary

\vspace{12pt}
\textbf{Corollary 3.2}.Let $g_{1},g_{2}$ be two Riemannian metrics on $M$,For any $\Omega \in \mathcal{C}(M)$ and any open subset $U\subset M$,we have $\mathbb{M}_U(\partial \Omega)_{g_{2}}\leq \sup_{x \in M}\frac{\sqrt{\det(g_{2})(x)}}{\sqrt{\det(g_{1})(x)}}  \mathbb{M}_U(\partial \Omega)_{g_{1}}.$

\vspace{12pt}
The proof of the next proposition follows closely from  Lemma 2.1 in [7].

\vspace{12pt}
\textbf{Proposition 3.3}.$W^{c}(M,g)$ depends continuously on the metric $g$(in the $C^{0}$ topology).
\begin{proof}
Suppose there is a sequence of Riemannian metrics \{$g_{i}$\}$_{i \in \mathbb{N}}$ that converges to $g$ in $C^{0}$ topology.Given $\epsilon >0$,let $\phi:[0,1]\to \mathcal{C}(M)$ is a sweepout such that
\begin{equation}
    \sup_{t \in [0,1]} \mathcal{A}^{c}_{g}(\phi(t))\leq W^{c}(M,g)+\epsilon.\tag{3.3}
\end{equation}
where $\mathcal{A}^{c}_{g}(\Omega)$ is the $\mathcal{A}^{c}$ of $\Omega$ with respect to $g$.

\vspace{5pt}
Since the flat norm of a Caccioppoli set $\Omega$ is equal to the $(n+1)$-dimensional Hausdorff measure of $\Omega$,then by Lemma 3.1 and Corollary 3.2 we see that the definition of sweepout is independent of choice of metrics.Thus by (3.1),(3.3)and Corollary 3.2 we have
\begin{align*}
W^{c}(M,g_{i})&\leq \sup_{t \in [0,1]} \mathcal{A}^{c}_{g_{i}}(\phi(t))\\
&\leq  \inf_{x \in M}\frac{\sqrt{\det(g_{i})(x)}}{\sqrt{\det(g)(x)}} \cdot \sup_{t \in [0,1]} \mathcal{A}^{c}_{g}(\phi(t))+(\sup_{x \in M}\frac{\sqrt{\det(g_{i})(x)}}{\sqrt{\det(g)(x)}}-\\  &\inf_{x \in M}\frac{\sqrt{\det(g_{i})(x)}}{\sqrt{\det(g)(x)}})\cdot(\sup_{t \in [0,1]} \mathbb{M}(\phi(t))_{g})\\
&\leq \inf_{x \in M}\frac{\sqrt{\det(g_{i})(x)}}{\sqrt{\det(g)(x)}} \cdot (W^{c}(M,g)+\epsilon)+(\sup_{x \in M}\frac{\sqrt{\det(g_{i})(x)}}{\sqrt{\det(g)(x)}}-\\  &\inf_{x \in M}\frac{\sqrt{\det(g_{i})(x)}}{\sqrt{\det(g)(x)}})\cdot(\sup_{t \in [0,1]} \mathbb{M}(\phi(t))_{g}). \tag{3.4}
\end{align*}
Since $\sup_{t \in [0,1]} \mathbb{M}(\phi(t))_{g}$ is finite by (3.3) and $\epsilon >0$ is arbitrary,we get 
\begin{equation}
\text{limsup}_{i\to \infty}W^{c}(M,g_{i})\leq W^{c}(M,g).\tag{3.5}
\end{equation}
Similarly,let $\phi_{i}:[0,1]\to \mathcal{C}(M)$ be a sequence of sweepout such that
\begin{equation}
      \sup_{t \in [0,1]} \mathcal{A}^{c}_{g_{i}}(\phi_{i}(t))\leq W^{c}(M,g_{i})+\epsilon.\tag{3.6}
\end{equation}
Then by the same reason we have
\begin{align*}
    W^{c}(M,g)&\leq \inf_{x \in M}\frac{\sqrt{\det(g)(x)}}{\sqrt{\det(g_{i})(x)}} \cdot (W^{c}(M,g_{i})+\epsilon)+(\sup_{x \in M}\frac{\sqrt{\det(g)(x)}}{\sqrt{\det(g_{i})(x)}}-\\  &\inf_{x \in M}\frac{\sqrt{\det(g)(x)}}{\sqrt{\det(g_{i})(x)}})\cdot(\sup_{t \in [0,1]} \mathbb{M}(\phi_{i}(t))_{g_{i}}). \tag{3.7}
\end{align*}
By (3.1),(3.5),(3.6) we have $\{\sup_{t \in [0,1]} \mathbb{M}(\phi_{i}(t))_{g_{i}}\}_{i \in \mathbb{N}}$ is bounded,since $\epsilon >0$ is arbitrary,we have 
\begin{equation}
     W^{c}(M,g)\leq \text{liminf}_{i\to \infty}W^{c}(M,g_{i}) \tag{3.8}
\end{equation}
\end{proof}

For any $c>0$,we denote by $\mathcal{P}^{c}_{g}$the class of smooth,closed,almost embedded $c-CMC$ hypersurface and $\widetilde{\mathcal{P}}^{c}_{g}$ represents the class of open set $\Omega$ such that $\partial \Omega$ is a smooth,closed,almost embedded $c-CMC$ hypersurface with respect to the inward unit normal of $\Omega$ in $(M,g)$.We also denote by $\overline{\mathcal{P}}^{c}_{g}\subset \widetilde{\mathcal{P}}^{c}_{g}$ the class of open set $\Omega$ such that $\mathcal{A}^{c}_{g}(\Omega)\leq W^{c}(M,g)$.

\vspace{12pt}
\textbf{Proposition 3.4}. Let $M^{n+1}(3\leq n+1\leq7)$ be a closed smooth manifold and let $g$ be a riemannian metric on $M$.Suppose $\#$ $\mathcal{P}^{c}_{g}<\infty$ and $\#$$\widetilde{\mathcal{P}}^{c}_{g}=N$,then there exists a Riemannian metric $h$ on $M$ such that 
\begin{enumerate}
\item[\textbf{(1)}] $\#$ $\widetilde{\mathcal{P}}^{c}_{h}\geq N+1$;
\item[\textbf{(2)}] $\mathcal{P}^{c}_{g}\subset$ $\mathcal{P}^{c}_{h}$;
\item[\textbf{(3)}] $||h-g||_{L^{\frac{n+1}{2}}(M,g)}\leq a(n) (\frac{W^{c}(M,g)}{c}N)^{\frac{2}{n+1}}$.
\end{enumerate}
where $a(n)$ is a constant depending only  on the dimension of $M$.

\vspace{12pt}
\begin{proof}
    Let $\{\Omega_{i}\}_{i=1,...,N}$ be the enumeration of $\widetilde{\mathcal{P}}^{c}_{g}$ and $\{\sum_{j}\}_{j=1,...,(N+m)}$  be the enumeration of $\mathcal{P}^{c}_{g}$.Note that $\partial \Omega_{k}\neq \partial\Omega_{l}$ if $k\neq l$.Since $\#$ $\widetilde{\mathcal{P}}^{c}_{g}$ and $\#$ $\mathcal{P}^{c}_{g}$ is finite,there exists connected open subset
$\{U_{i}\}_{i=1,...,N}$ in $M$ such that 
\begin{enumerate}
\item[\textbf{i)}] $U_{i}\subset \Omega_{i}$;
\item[\textbf{ii)}] $U_{k}\bigcap U_{l}=\varnothing$ for any $k\neq l$;
\item[\textbf{iii)}] The intersection of $U:=\bigcup_{i=1}^{N} U_{i}$ and $\bigcup_{i=1}^{N+m}\sum_{i}$ is empty.
\end{enumerate}
Let $\{f_{i}\}_{i=1,...,N}$ be the smooth function on $M$ such that
\begin{enumerate}
\item[\textbf{i)}] supp$f_{i}\subset U_{i}$;
\item[\textbf{ii)}] $f_{i}(x)\geq 0$ for any $x\in M$;
\item[\textbf{iii)}] $\int_{U_{i}}|f_{i}|^{\frac{n+1}{2}}=1$.
\end{enumerate}
Define $f=\sum_{i=1}^{N}f_{i}$ ,$g_{t} = (1+tf) g$  for $t\geq 0$.By the construction of $\{f_{i}\}_{i=1,...,N}$,we have $\mathcal{P}^{c}_{g}\subset \mathcal{P}^{c}_{g_{t}}$.The following argument is similar to the proof of proposition 3.1 in [7].Since $g|_{M\setminus U}=g_{t}|_{M\setminus U}$,then we have $\mathbb{M}_{g}(\partial\Omega_{i})=\mathbb{M}_{g_{t}}(\partial\Omega_{i})$ for every $\Omega_{i} \in \widetilde{\mathcal{P}}^{c}_{g}$.Moreover,we have
\begin{align*}
    \mathcal{A}^{c}_{g_{t}}(\Omega_{i})&=\mathcal{A}^{c}_{g}(\Omega_{i})-c(\int_{\Omega_{i}}(1+tf)^{\frac{n+1}{2}}d\mathcal{H}^{n+1}_{g}-\int_{\Omega_{i}}1 d\mathcal{H}^{n+1}_{g})\\
    &\leq \mathcal{A}^{c}_{g}(\Omega_{i})-c\int_{\Omega_{i}}(tf)^{\frac{n+1}{2}}d\mathcal{H}^{n+1}_{g}\\
    &\leq \mathcal{A}^{c}_{g}(\Omega_{i})-ct^{\frac{n+1}{2}}.\tag{3.9}
\end{align*}
In the following, we show that there exists $\widetilde{t}>0$ such that $\widetilde{\mathcal{P}}_{g_{\widetilde{t}}}^{c}$ contains some $\widetilde{\Omega}$ with $\partial\widetilde{\Omega}\cap U\neq\emptyset$ and the $L^{\frac{n+1}{2}}$ distance between $g$ and $g_{\widetilde{t}}$ possesses an upper bound as stated in item (3) of the proposition.Let $t_{0}>0$,suppose for any $t\in[0,t_{0}]$,every $\Omega \in \widetilde{\mathcal{P}}^{c}_{g_{t}}$ satisfies that $\partial \Omega \subset M \setminus U$.We will show 
\begin{equation}
    t_{0}< (\frac{W^{c}(M,g)}{c})^{\frac{2}{n+1}}.\tag{3.10}
\end{equation}

By our assumption on $t_{0}$ and Theorem 2.7,for all $t\in [0,t_{0}]$ we have  $W^{c}(M,g_{t})=\mathcal{A}^{c}_{g_{t}}(\Omega_{i_{t}})$ for some $\Omega_{i_{t}}\in\widetilde{\mathcal{P}}^{c}_{g}$.We will prove (3.10) in the following two cases respectively.

\vspace{8pt}
\textbf{\textbf{Case 1:}} Suppose there exits $\Omega_{i_{t_{0}}}\in \overline{\mathcal{P}}^{c}_{g}$ such that $W^{c}(M,g_{t_{0}})=\mathcal{A}^{c}_{g_{t_{0}}}(\Omega_{i_{t_{0}}})$.
Then we know from Proposition 2.5 that $\mathcal{A}^{c}_{g_{t_{0}}}(\Omega_{i_{t_{0}}})>0$.Hence (3.9) implies (3.10) is true.

\vspace{8pt}
 \textbf{\textbf{Case 2:}} Suppose $W^{c}(M,g_{t_{0}})$ only can be realized by the $\mathcal{A}^{c}_{g_{t_{0}}}(\Omega_{i_{t_{0}}})$ for some $\Omega_{i_{t_{0}}}\in  \widetilde{\mathcal{P}}^{c}_{g} \setminus \overline{\mathcal{P}}^{c}_{g}$,i.e.
 \begin{equation*}
     \mathcal{A}^{c}_{g}(\Omega_{i_{t_{0}}})> W^{c}(M,g).
 \end{equation*}
 Then we claim that the following statements are true:

 \vspace{8pt}

 \item[\textbf{Claim 3.5}]: There exists $\{t_{j}\}_{j=1,...,k}\subset [0,t_{0}]$ with $t_{1}\leq ...\leq t_{k}$ and exists $\{\Omega_{j}\}_{j=1...k}\in \widetilde{\mathcal{P}}^{c}_{g}$ satisfies
 \begin{enumerate}
     \item [(i)] $\Omega_{1}\in \overline{\mathcal{P}}^{c}_{g},\Omega_{2}\in\widetilde{\mathcal{P}}^{c}_{g} \setminus \overline{\mathcal{P}}^{c}_{g} $;
     \item [(ii)] $\mathcal{A}^{c}_{g_{t_{j}}}(\Omega_{j})=\mathcal{A}^{c}_{g_{t_{j}}}(\Omega_{j+1})=W^{c}(M,g_{t_{j}})$ for $j=1,...,(k-1)$;
     \item [(iii)] $\mathcal{A}^{c}_{g_{t_{k}}}(\Omega_{k})=\mathcal{A}^{c}_{g_{t_{k}}}(\Omega_{i_{t_{0}}})=W^{c}(M,g_{t_{k}})$.
 \end{enumerate}
 
\vspace{8pt}
 Define $\sigma:[0,t_{0}]\to \mathbb{R}^{2}$ by 
  \begin{equation*}
     \sigma(t)=(t,W^{c}(M,g_{t}))
 \end{equation*}
 Define $\psi_{i}:[0,t_{0}]\to \mathbb{R}^{2},i=1,...,N$ by
 \begin{equation*}
     \psi_{i}(t)=(t,\mathcal{A}^{c}_{g_{t}}(\Omega_{i}))
 \end{equation*}
 Denotes the image of $\psi_{i}$ by $\Gamma_{i}$.Intuitively,claim 3.5 describes the evolutionary process of the image of $\sigma$ among $\{\Gamma_{i}\}_{i=1,...,N}$ as $t$ increases,which  starts at  $(0,W^{c}(M,g))$ and ends at $(t_{0},\mathcal{A}^{c}_{g_{t_{0}}}(\Omega_{i_{t_{0}}}))$.

\vspace{10pt}
From Proposition 3.3 and the assumption on  $t_{0}$,we deduce that the image of $\sigma$  is entirely  contained in a path-connected component of $\bigcup_{i=1}^{N}\Gamma_{i}$.More precisely,the image of $\sigma$ is contained in the path-connected component $C$ that contains some $\Gamma_{i_{0}}$ such that
 \begin{equation*}
     \mathcal{A}^{c}_{g}(\Omega_{i_{0}})=W^{c}(M,g)
 \end{equation*}
By the fact that $W^{c}(M,g_{t_{0}})=\mathcal{A}^{c}_{g_{t_{0}}}(\Omega_{i_{t_{0}}})$,together with the continuity of $W^{c}(M,g_{t})$,we deduce that $\Gamma_{i_{t_{0}}} \subset C$  and Claim 3.5 is true.

\vspace{8pt}
Throughout the following proof,we will use the notations used in Claim 3.5.Firstly for every $i\in\{1,...,N\}$ and $t\in[0,t_{0}]$,we have
\begin{align*}
\frac{d \mathcal{A}^{c}_{g_{t}}(\Omega_{i})}{dt}&=-c\int_{\Omega_{i}}\frac{n+1}{2}(1+tf)^{\frac{n-1}{2}}f d\mathcal{H}^{n+1}_{g}\\  
&\leq -c\int_{\Omega_{i}}\frac{n+1}{2}t^{\frac{n-1}{2}}f^{\frac{n+1}{2}}d\mathcal{H}^{n+1}_{g}\\ 
&= -c\frac{n+1}{2}t^{\frac{n-1}{2}}\\
&=\frac{d(-ct^{\frac{n+1}{2}})}{dt}\tag{3.11}
\end{align*}
Then by (3.9) and Claim 3.5,
\begin{align*}
\mathcal{A}^{c}_{g_{t_{1}}}(\Omega_{2})=\mathcal{A}^{c}_{g_{t_{1}}}(\Omega_{1})\leq W^{c}(M,g)-ct_{1}^{\frac{n+1}{2}} \tag{3.12}
\end{align*}
Applying (3.11) we obtain
\begin{align*}
\mathcal{A}^{c}_{g_{t_{2}}}(\Omega_{3})=\mathcal{A}^{c}_{g_{t_{2}}}(\Omega_{2})\leq  W^{c}(M,g)-ct_{2}^{\frac{n+1}{2}}  \tag{3.13}  
\end{align*}
By an induction argument we deduce that
\begin{align*}
\mathcal{A}^{c}_{g_{t_{k}}}(\Omega_{i_{t_{0}}})=\mathcal{A}^{c}_{g_{t_{k}}}(\Omega_{k})\leq  W^{c}(M,g)-ct_{k}^{\frac{n+1}{2}}  \tag{3.14}  
\end{align*}
Finally  by (3.11) we have
\begin{equation*}
 W^{c}(M,g_{t_{0}})=\mathcal{A}^{c}_{g_{t_{o}}}(\Omega_{i_{t_{0}}})\leq  W^{c}(M,g)-ct_{0}^{\frac{n+1}{2}}   \tag{3.15}
\end{equation*}
Combining (3.15) and Proposition 2.5 we deduce that (3.10) holds.Thus there must exist $\widetilde{t}\leq (\frac{W^{c}(M,g)}{c})^{\frac{2}{n+1}}$ such that $\widetilde{\mathcal{P}}_{g_{\widetilde{t}}}^{c}$ contains some $\widetilde{\Omega}$ with $\partial\widetilde{\Omega}\cap U\neq\emptyset$.Let $h:=g_{\widetilde{t}}$,by the fact $g|_{M\setminus U}=g_{t}|_{M\setminus U}$ for every $t>0$ and  $U\bigcap_{i=1}^{N+m} \sum_{i}=\emptyset$,we have items (1) and (2) of the proposition hold.To show part (3) of the Proposition,we first compute the $L^{\frac{n+1}{2}}$ distance between $g$ and $g_{(\frac{W^{c}(M,g)}{c})^{\frac{2}{n+1}}}$.By the definition of $g_{(\frac{W^{c}(M,g)}{c})^{\frac{2}{n+1}}}$,we have
\begin{align*}
 ||g_{(\frac{W^{c}(M,g)}{c})^{\frac{2}{n+1}}}-g||_{L^{\frac{n+1}{2}}(M,g)}&=(\int_{M}|g_{(\frac{W^{c}(M,g)}{c})^{\frac{2}{n+1}}}-g|_{g}^{\frac{n+1}{2}}d\mathcal{H}^{n+1}_{g})^{\frac{2}{n+1}} \\
 &=(\frac{W^{c}(M,g)}{c})^{\frac{2}{n+1}}(\int_{M}|fg|_{g}^{\frac{n+1}{2}}d\mathcal{H}^{n+1}_{g})^{\frac{2}{n+1}}\\
 &=a(n) (\frac{W^{c}(M,g)}{c}N)^{\frac{2}{n+1}}
\end{align*}

Since the $L^{\frac{n+1}{2}}$ distance between $g$ and $g_{t}$ is increasing in $t$ and $\widetilde{t}\leq (\frac{W^{c}(M,g)}{c})^{\frac{2}{n+1}}$,conclusion (3) of Proposition 3.4 follows.
\end{proof}
\vspace{15pt}

$Proof$ $of$ $Theorem$ $1.1$:

\vspace{15pt}
Since Proposition 3.4 holds,we only need to deal with the case when $k\geq 2$.By the the proof of Proposition 3.4,there exists a smooth and nonnegative function $f_{1}$ on $M$,such that the intersection of the union of all $c-CMC$ hypersurface in $(M,g)$ with supp$f_{1}$ is empty and there exits $t\in[0,(\frac{W^{c}(M,g)}{c})^{\frac{2}{n+1}}]$ with  $\#\widetilde{\mathcal{P}}^{c}_{(1+tf_{1}).g}\geq  N+1$.Define $\widetilde{t_{1}}$ as 
\begin{equation*}
 \widetilde{t_{1}}=\text{inf}\{t | \#\widetilde{\mathcal{P}}^{c}_{(1+tf_{1}).g}\geq  N+1\}    
\end{equation*}
\textbf{Case 1:} Suppose that $\widetilde{t_{1}}>0$,then by the construction of $f_{1}$ there exists $t_{1}>0$ and $\epsilon_{1} >0$ such that if we denote $(1+t_{1}f_{1}).g$ as $g_{1}$ then we have
\begin{enumerate}[itemsep=0.3em, parsep=0.5em] 
    \item [(1)] $\mathcal{P}^{c}_{g}\subset$ $\mathcal{P}^{c}_{g_{1}}$;
    \item [(2)] $\#\widetilde{\mathcal{P}}^{c}_{g_{1}}\geq  N+1$ ;
    \item [(3)] $W^{c}(M,g)-W^{c}(M,(1+t_{1}f_{1}).g)=\epsilon_{1}$;
    \item [(4)] $||g_{1}-g||_{L^{\frac{n+1}{2}}(M,g)}\leq a(n) (\frac{W^{c}(M,g)}{c}N)^{\frac{2}{n+1}}$.
\end{enumerate}
It follows from the process above that there exists positive real numbers $\{t_{i}\}_{i=2....,n(k)}$ and smooth nonnegative functions $\{f_{i}\}_{i=2,...,n(k)}$ on $M$,where $n(k)\leq k$,such that if we denote $g_{i}=\prod_{j=2}^{j=i}(1+t_{i}f_{i}).g_{1}$  for every $i\geq 2$,then we have
\begin{enumerate}[itemsep=0.3em, parsep=0.5em] 
    \item [(a)] $\mathcal{P}^{c}_{g_{i-1}}\subset$ $\mathcal{P}^{c}_{g_{i}}$;
    \item [(b)] $W^{c}(M,g_{i})-W^{c}(M,g_{i-1})\leq \frac{\epsilon_{1}}{k}$;
    \item [(c)] $||g_{i}-g_{i-1}||_{L^{\frac{n+1}{2}}(M,g_{i-1})}\leq a(n) [\frac{W^{c}(M,g_{i-1})}{c}(\#\widetilde{\mathcal{P}}^{c}_{g_{i-1}})]^{\frac{2}{n+1}}$;
     \item [(d)] $\#\widetilde{\mathcal{P}}^{c}_{g_{i}}\geq \#\widetilde{\mathcal{P}}^{c}_{g_{i-1}} +1$,$\#\widetilde{ \mathcal{P}}^{c}_{g_{n(k)-1}}< \#\widetilde{\mathcal{P}}^{c}_{g} +k$,$\#\widetilde{ \mathcal{P}}^{c}_{g_{n(k)}}\geq \#\widetilde{\mathcal{P}}^{c}_{g} +k$.
\end{enumerate}
Since 
\begin{align*}
||g_{i}-g_{i-1}||_{L^{\frac{n+1}{2}}(M,g_{i-1})}&=(\int_{M}|g_{i}-g_{i-1}|_{g_{i-1}}^{\frac{n+1}{2}}d\mathcal{H}^{n+1}_{g_{i-1}})^{\frac{2}{n+1}} \notag \\  
&=(\int_{M}[\prod_{j=1}^{i-1}(1+t_{j}f_{j})t_{i}f_{i}|g|_{g_{i-1}}]^{\frac{n+1}{2}}d\mathcal{H}^{n+1}_{g_{i-1}})^{\frac{2}{n+1}} \notag \\  
&=(\int_{M}[\prod_{j=1}^{i-1}(1+t_{j}f_{j})t_{i}f_{i}|g|_{g}]^{\frac{n+1}{2}}[\prod_{j=1}^{i-1}(1+t_{j}f_{j})]^{-\frac{n+1}{2}}d\mathcal{H}^{n+1}_{g_{i-1}})^{\frac{2}{n+1}} \notag \\  
&=(\int_{M}[\prod_{j=1}^{i-1}(1+t_{j}f_{j})t_{i}f_{i}|g|_{g}]^{\frac{n+1}{2}}d\mathcal{H}^{n+1}_{g})^{\frac{2}{n+1}} \notag \\  
&=(\int_{M}|g_{i}-g_{i-1}|_{g}^{\frac{n+1}{2}}d\mathcal{H}^{n+1}_{g})^{\frac{2}{n+1}} \notag \\  
&=||g_{i}-g_{i-1}||_{L^{\frac{n+1}{2}}(M,g)}. \tag{3.16}  
\end{align*}
Thus we have
\begin{align*}
  ||g_{n(k)}-g||_{L^{\frac{n+1}{2}}(M,g)}&\leq \sum_{i=2}^{n(k)}||g_{i}-g_{i-1}||_{L^{\frac{n+1}{2}}(M,g)}+||g_{1}-g||_{L^{\frac{n+1}{2}(M,g)}}\\
  &=\sum_{i=2}^{n(k)}||g_{i}-g_{i-1}||_{L^{\frac{n+1}{2}}(M,g_{i-1})}+||g_{1}-g||_{L^{\frac{n+1}{2}(M,g)}}\\
  &\leq \sum_{i=2}^{n(k)}a(n) [\frac{W^{c}(M,g_{i-1})}{c}(\#\widetilde{\mathcal{P}}^{c}_{g_{i-1}})]^{\frac{2}{n+1}}+a(n) (\frac{W^{c}(M,g)}{c}N)^{\frac{2}{n+1}}\\
  &\leq \sum_{i=0}^{k-1}a(n) [\frac{W^{c}(M,g)}{c}(N+i)]^{\frac{2}{n+1}}\\
  \text{(by h\"{o}lder inequality)}&\leq a(n)[\frac{W^{c}(M,g)}{c}]^{\frac{2}{n+1}}k(N+\frac{k-1}{2})^{\frac{2}{n+1}}. \tag{3.17}
\end{align*}
\textbf{Case 2:} Suppose that $\widetilde{t_{1}}=0$.In this case,$W^{c}(M,(1+tf_{1}).g)$ might increase as t increases and thus (3) in Case 1 may not hold anymore.Therefore, some modifications need to be made here compared with Case 1.By  induction and (3.16) we have one of the following two statements is true. 
\begin{enumerate}

\item [(i)] There exist $\{t_{i}\}_{i=1,...,m(k)}>0$,$\{\epsilon_{i}\}_{i=1,...,m(k)}>0$ and smooth nonnegative functions $\{f_{i}\}_{i=1,...,m(k)}$ on $M$.For each $i$,let $g_{i}=\prod_{j=1}^{j=i}(1+t_{i}f_{i}).g_{1}$  with $g_{0}=g$,then we have
\begin{enumerate}[itemsep=0.5em, parsep=0.5em] 
   
    \item [(1)] $\mathcal{P}^{c}_{g_{i-1}}\subset$ $\mathcal{P}^{c}_{g_{i}}$;
    \item [(2)] $\#\widetilde{\mathcal{P}}^{c}_{g_{i}}\geq  \#\widetilde{\mathcal{P}}^{c}_{g_{i-1}}+1$,$\#\widetilde{ \mathcal{P}}^{c}_{g_{m(k)-1}}< \#\widetilde{\mathcal{P}}^{c}_{g} +k$,$\#\widetilde{ \mathcal{P}}^{c}_{g_{m(k)}}\geq \#\widetilde{\mathcal{P}}^{c}_{g} +k$;
     \item [(3)] $ \widetilde{t_{i}}:=\text{inf}\{t | \#\widetilde{\mathcal{P}}^{c}_{(1+tf_{i}).g_{i-1}}\geq  \#\widetilde{\mathcal{P}}^{c}_{g_{i-1}}+1\}=0    $ for every $i$;
     \item [(4)] $||g_{i}-g_{i-1}||_{L^{\frac{n+1}{2}}(M,g)}=\epsilon_{i}$;
    \item [(5)] $\sum_{i=1}^{m(k)}\epsilon_{i}\leq  a(n)[\frac{W^{c}(M,g)}{c}]^{\frac{2}{n+1}}k(N+\frac{k-1}{2})^{\frac{2}{n+1}}. $.
    
\end{enumerate}

\vspace{5pt}
\item[(ii)]There exist $\{t_{i}\}_{i=1,...,m(k)}>0$,$\{\epsilon_{i}\}_{i=1,...,m(k)-m}>0$(where $m>0$) and smooth nonnegative functions $\{f_{i}\}_{i=1,...,m(k)}$ on $M$.For each $i$,let $g_{i}=\prod_{j=1}^{j=i}(1+t_{i}f_{i}).g_{1}$  with $g_{0}=g$,then we have
\begin{enumerate}[itemsep=0.5em, parsep=0.5em] 
   
    \item [(1)] $\mathcal{P}^{c}_{g_{i-1}}\subset$ $\mathcal{P}^{c}_{g_{i}}$;
    \item [(2)] $\#\widetilde{\mathcal{P}}^{c}_{g_{i}}\geq  \#\widetilde{\mathcal{P}}^{c}_{g_{i-1}}+1$ , $\#\widetilde{ \mathcal{P}}^{c}_{g_{m(k)-1}}< \#\widetilde{\mathcal{P}}^{c}_{g} +k$ , $\#\widetilde{ \mathcal{P}}^{c}_{g_{m(k)}}\geq \#\widetilde{\mathcal{P}}^{c}_{g} +k$ ;
     \item [(3)] $ \widetilde{t_{i}}:=\text{inf}\{t | \#\widetilde{\mathcal{P}}^{c}_{(1+tf_{i}).g_{i-1}}\geq  \#\widetilde{\mathcal{P}}^{c}_{g_{i-1}}+1\}=0    $ for $i\leq m(k)-m$,and $ \widetilde{t_{i}}:=\text{inf}\{t | \#\widetilde{\mathcal{P}}^{c}_{(1+tf_{i}).g_{i-1}}\geq  \#\widetilde{\mathcal{P}}^{c}_{g_{i-1}}+1\}>0    $ for $i=m(k)-m+1$;
     \item [(4)] $||g_{i}-g_{i-1}||_{L^{\frac{n+1}{2}}(M,g)}=\epsilon_{i}$ for $i=1,...,(m(k)-m)$,$||g_{i}-g_{i-1}||_{L^{\frac{n+1}{2}}(M,g)}\leq a(n) [\frac{W^{c}(M,g_{m(k)-m})}{c}(\#\widetilde{\mathcal{P}}^{c}_{g_{i-1}})]^{\frac{2}{n+1}}$ for $i=(m(k)-m+1),...,k$
    
    \item [(5)] $a(n) [\frac{W^{c}(M,g_{i-1})}{c}(N+i-1)]^{\frac{2}{n+1}}-\epsilon_{i}\geq \sum_{j=i}^{k-1}a(n)\bigl[[\frac{W^{c}(M,g_{i})}{c}(N+j)]^{\frac{2}{n+1}}- [\frac{W^{c}(M,g_{i-1})}{c}(N+j)]^{\frac{2}{n+1}} \bigr]$ for $i=1,...,(m(k)-m)$;
   
\end{enumerate}

\end{enumerate}

\vspace{10pt}
Firstly suppose (i) holds,then

\begin{align*}
 ||g_{m(k)}-g_{0}||_{L^{\frac{n+1}{2}}(M,g)}&\leq \sum_{i=1}^{m(k)}||g_{i}-g_{i-1}||_{L^{\frac{n+1}{2}}(M,g)}\\
 &=\sum_{i=1}^{m(k)}\epsilon_{i}\leq a(n)[\frac{W^{c}(M,g)}{c}]^{\frac{2}{n+1}}k(N+\frac{k-1}{2})^{\frac{2}{n+1}}.
\end{align*}

Now suppose (ii) holds,then we have

\begin{align*}
 ||g_{m(k)}-g_{0}||_{L^{\frac{n+1}{2}}(M,g)}&\leq \sum_{i=1}^{m(k)-m}||g_{i}-g_{i-1}||_{L^{\frac{n+1}{2}}(M,g)}+\sum_{i=m(k)-m+1}^{m(k)}||g_{i}-g_{i-1}||_{L^{\frac{n+1}{2}}(M,g)}\\
 &\leq  \sum_{i=1}^{m(k)-m}\epsilon_{i}+\sum_{i=m(k)-m+1}^{m(k)}a(n) [\frac{W^{c}(M,g_{m(k)-m})}{c}(\#\widetilde{\mathcal{P}}^{c}_{g_{i-1}})]^{\frac{2}{n+1}}\\
& \leq \sum_{i=1}^{m(k)-m}\bigl[\sum_{j=i-1}^{k-1}a(n)[\frac{W^{c}(M,g_{i-1})}{c}(N+j)]^{\frac{2}{n+1}}-[\sum_{j=i}^{k-1}a(n)[\frac{W^{c}(M,g_{i})}{c}(N+j)]^{\frac{2}{n+1}}\bigr]\\
&+\sum_{i=m(k)-m+1}^{m(k)}a(n) [\frac{W^{c}(M,g_{m(k)-m})}{c}(\#\widetilde{\mathcal{P}}^{c}_{g_{i-1}})]^{\frac{2}{n+1}}\\
&=\sum_{i=0}^{k-1}a(n)[\frac{W^{c}(M,g)}{c}(N+i)]^{\frac{2}{n+1}}-\sum_{i=m(k)-m}^{k-1}a(n)[\frac{W^{c}(M,g_{m(k)-m})}{c}(N+i)]^{\frac{2}{n+1}}+\\
&+\sum_{i=m(k)-m+1}^{m(k)}a(n) [\frac{W^{c}(M,g_{m(k)-m})}{c}(\#\widetilde{\mathcal{P}}^{c}_{g_{i-1}})]^{\frac{2}{n+1}}\\
&\leq \sum_{i=0}^{k-1}a(n)[\frac{W^{c}(M,g)}{c}(N+i)]^{\frac{2}{n+1}}\\
\text{(by h\"{o}lder inequality)}&\leq a(n)[\frac{W^{c}(M,g)}{c}]^{\frac{2}{n+1}}k(N+\frac{k-1}{2})^{\frac{2}{n+1}}. 
\end{align*}
\qed

        \begin{flushleft}
			\medskip\noindent
			\begin{tabbing}
				XXXXXXXXXXXXXXXXXXXXXXXXXX*\=\kill
				Xiaoxiang Jiao\\
				School of Mathematical Sciences, University of Chinese Academy of Sciences\\
				19A Yuquan Road, Beijing, 100049, China\\
				
				E-mail: xxjiao@ucas.ac.cn
				
			\end{tabbing}
		\end{flushleft}

        \begin{flushleft}
			\medskip\noindent
			\begin{tabbing}
				XXXXXXXXXXXXXXXXXXXXXXXXXX*\=\kill
				Wenduo Zou\\
				School of Mathematical Sciences, University of Chinese Academy of Sciences\\
				19A Yuquan Road, Beijing, 100049, China\\
				
				E-mail: zouwenduo23@mails.ucas.ac.cn
				
			\end{tabbing}
		\end{flushleft}

\medskip \medskip \medskip
		\noindent
		\vskip 0.3in
		\begin{thebibliography}{99}
			\vskip 0.2in
			
			\bibitem[1]{1}
			F. Almgren.
			\newblock {The theory of varifolds}.
			\newblock {\em Mimeographed notes,Princeton}, 1965.
			
			\bibitem[2]{2}
			V. I. Arnold.
			\newblock {\em Arnold’s problems}.
			\newblock {Springer, 2004}.
			
			\bibitem[3]{3}
			A. Dey.
			\newblock {Existence of multiple closed CMC hypersurfaces with small mean curvature}.
			\newblock {\em J. Differential Geom}~,125(2), 379–403, 2023.
			
			\bibitem[4]{4}
			P. Glynn-Adey, Y. Liokumovich. 
			\newblock {Ricci curvature, and minimal hypersurfaces}.
			\newblock {\em J. Differential Geom},105(1), 33–54, 2017.
			
			\bibitem[5]{5}
			M. Gromov.In
			\newblock {Dimension, nonlinear spectra and width}.In
			\newblock {\em Geometric aspects of functional analysis (1986/87)},volume 1317 of \newblock{\em Lecture Notes in Math.},1317, 132–184, Springer, Berlin, 1988.


            \bibitem[6]{6}
			M. Gromov.
			\newblock {Isoperimetry of waists and concentration of maps}.
			\newblock {\em Geom. Funct. Anal}, 13(1), 178–215, 2003.
			
			\bibitem[7]{7}
			K. Irie, F. C. Marques, A. Neves. 
			\newblock {Density of minimal hypersurfaces for generic metrics}.
			\newblock {\em Ann. of Math},963–972, 2018
			
			
			\bibitem[8]{8}
			Y. Li. 
			\newblock{ Existence of infinitely many minimal hypersurfaces in higher-dimensional closed manifolds with generic metrics}.
            \newblock{\em J. Differential Geom},124(2), 381–395, 2023

            
			\bibitem[9]{9}
			Y. Liokumovich, F. C. Marques, A. Neves.
			\newblock {Weyl law for the volume spectrum}.
			\newblock {\em Ann. of Math}, 187(3), 933–961, 2018.
			
			\bibitem[10]{10}
			L. Mazurowski, X. Zhou.
			\newblock {The half-volume spectrum of a manifold}.
			\newblock {\em   Calc. Var. Partial Differential Equations}, 64, 155, 2025.
			
			\bibitem[11]{11}
			L. Mazurowski, X. Zhou. 
			\newblock {Infinitely many half-volume constant mean curvature hypersurfaces via min-max theory}.
			\newblock {\em arXiv preprint arXiv:2405.00595,}2024.
			
			\bibitem[12]{12}
			L. Mazurowski, X. Zhou.
			\newblock {An alternative for constant mean curvature hypersurfaces}.
			\newblock {\em arXiv preprint arXiv:2408.13864,}2024.
			
			\bibitem[13]{13}
			R. Mazzeo, F. Pacard.
			\newblock {Foliations by constant mean curvature tubes}.
			\newblock {\em Communications in analysis and geometry}, 13(4), 633–670, 2005.
			
			\bibitem[14]{14}
			F. Pacard and X. Xu.
			\newblock {Constant mean curvature spheres in riemannian manifolds}.
			\newblock {\em Manuscripta mathematica}, 128:275--295, 2009.
			
			\bibitem[15]{15}
			J. Pitts.
			\newblock {\em Existence and regularity of minimal surfaces on Riemannian manifolds}.
			\newblock Princeton University Press, 1981.
			
			\bibitem[16]{16}
			R. Schoen and L. Simon.
			\newblock {Regularity of stable minimal hypersurfaces}.
			\newblock {\em Comm. Pure Appl. Math.}, 34:741--797, 1981.
			
			\bibitem[17]{17}
			L. Simon.
			\newblock {\em Lectures on geometric measure theory}.
			\newblock {\em Proceedings of the Centre for Mathematical Analysis}, Australian National University, 3, 1983.
			
			\bibitem[18]{18}
			A. Song.
			\newblock {Existence of infinitely many minimal hypersurfaces in closed manifolds}.
			\newblock {\em Ann. of Math.}, 197(3):859--895, 2023.
			
			\bibitem[19]{19}
			Z. Wang.
			\newblock {Conformal upper bounds for the volume spectrum}.
			\newblock {\em Geom. Funct. Anal.}, 31:992--1012, 2021.
			
			\bibitem[20]{20}
			S.-T. Yau.Problem section.In
			\newblock {\em Seminar on differential geometry},
			volume 102 of  \newblock {\em Annals of Mathematics Studies}, Princeton University Press, 1982.
			
			\bibitem[21]{21}
			R. Ye.
			\newblock {Foliation by constant mean curvature spheres}.
			\newblock {\em Pacific Journal of Mathematics}, 147(2):381--396, 1991.
			
			\bibitem[22]{22}
			X. Zhou.
			\newblock {Min-max minimal hypersurface in $(M^{n+1},g)$ with $Ric> 0$ and $2\leq n \leq 6$}.
			\newblock {\em J. Differential Geom.}, 100(1):129--160, 2015.
			
			\bibitem[23]{23}
			X. Zhou.
			\newblock {Min-max hypersurface in manifold of positive Ricci curvature}.
			\newblock {\em J. Differential Geom.}, 105(2):291--343, 2017.
			
			\bibitem[24]{24}
			X. Zhou.
			\newblock {On the multiplicity one conjecture in min-max theory}.
			\newblock {\em Annals of Mathematics}, 192(3):767--820, 2020.
			
			\bibitem[25]{25}
			X. Zhou.
			\newblock {Mean curvature and variational theory}.In
			\newblock {\em ICM—International Congress of Mathematicians.Vol. IV. Sections 5–8}, 2696–2717, EMS Press, Berlin, 2023.
			
			\bibitem[26]{26}
			X. Zhou and J. Zhu.
			\newblock {Min-max theory for constant mean curvature hypersurfaces}.
			\newblock {\em Invent. Math.}, 218(2):441--490, 2019.
			
		\end{thebibliography}
\end{document}